\documentclass[12pt,twoside]{article}
\usepackage{amsmath,amssymb}
\usepackage{graphicx}

\newtheorem{lm}{Lemma}[section]
\newtheorem{thm}{Theorem}[section]

\numberwithin{equation}{section}

\newcounter{saveeqn}%

\makeatletter
\evensidemargin
\oddsidemargin
\makeatother

\title{\Large\bf
Smooth solutions of a class of iterative functional equations
\thanks{
Supported by NSFC $\#$12001537, the Start-up Funding of Chongqing Normal University 20XLB033, and the Research Project of Chongqing Education Commission CXQT21014.
}
}

\author{Weiwei Shi,~~Xiao Tang
\footnote
{Corresponding author. E-mail address: x.tang@cqnu.edu.cn (X. Tang).
}
	\\
{\small School of Mathematical Sciences, Chongqing Normal University,}
	\\
{\small Chongqing 401331, PR China}
}

\date{}

\begin{document}
\maketitle

\begin{abstract}
Imposing some conditions on derivatives of the known functions, using the Fiber Contraction Theorem we prove the existence of $C^1$ solutions of a class of iterative functional equations which involves iterates of the unknown functions and a nonlinear term. 

\vskip 0.2cm


\vskip 0.2cm
{\bf AMS (2020) subject classification:} 39B12; 39B22
\end{abstract}

\baselineskip 15pt   
\parskip 10pt         

\thispagestyle{empty}
\setcounter{page}{1}


\setcounter{equation}{0}
\setcounter{lm}{0}
\setcounter{thm}{0}
\setcounter{rmk}{0}
\setcounter{df}{0}
\setcounter{cor}{0}
\setcounter{exa}{0}
\allowdisplaybreaks[2]


\section{Introduction}
Iterative functional equations (\cite{BJ,KCG} and references therein) involving the iterates of unknown functions are extensively studied. In particular, iterative roots (\cite{LJLZ,LZ} and references therein) and polynomial-like iterative equations (\cite{MR,ZNX} and references therein) are of this class. The complication comes from the fact that the iteration operator is nonlinear.

Considering  the muli-variable functional equation
\begin{equation}
	x+\phi(y+\phi(x))=y+\phi(x+\phi(y)),
	\label{two-variable}
\end{equation}
N. Brillou\"{e}t-Belluot (\cite{Belluot}) in 2000 proposed the
second order iterative functional equation
\begin{equation}\label{sw}
	\phi^2(x)=\phi(x+a)-x
\end{equation}
in the problem session of the 38th ISFE held in Hungary, which was mentioned again by K. Baron (\cite{Baron}) in 2003.
It is easy to see that \eqref{two-variable} with $y=0$ is reduced to \eqref{sw} with $a=\phi(0)$.
Because it should be re-considered that the existence of solutions of functional equations even when they vary slightly,  a large number of researchers focus both on \eqref{two-variable} and \eqref{sw}; see 
\cite{Balcerowski,Yarczyk,Sablik,ZWN,ZYY,TX}. By \cite[Corollary 3.8]{DM} or \cite[Theorem 11]{Matkowski} or \cite[Theorem 5]{YDL}, the equation \eqref{sw} has no continuous solutions on $\mathbb{R}$ when $a=0$.

In 2010,  N. Brillou\"et-Belluot and W. Zhang (\cite{ZWN}) considered  a general one of the equation \eqref{sw},
\begin{equation}
	\phi^2(x)=\lambda \phi(x+a)+\mu x,
\label{general-1}
\end{equation}
where $\lambda,a$ and $\mu $ are real such that $a\lambda \neq 0$. They proved that there are Lipschitz solutions on any given compact interval of \eqref{general-1} under the condition 
\begin{equation*}
|\lambda|> \max \{2,2\sqrt{2|\mu|}\}~~~\text{and}~~~1+2|\mu|< |\lambda| \leq 2.
\end{equation*}
Moreover, piecewise continuous solutions on a bounded interval of \eqref{general-1} are constructed in the case that 
\begin{equation*}
0 \leq \mu < 1~~~~\text{and}~~~~ \lambda \geq 2(1-\mu) .
\end{equation*}
Besides, Y. Zeng and W. Zhang (\cite{ZYY}) showed the nonexistence of continuous solutions on $\mathbb{R}$ of \eqref{general-1} when $\lambda=1$ and $\mu\le -1$, which answers that the equation \eqref{sw} also has no continuous solutions on $\mathbb{R}$ in the case  $a\neq 0$, and proved the existence of continuous solutions on $\mathbb{R}$ of \eqref{general-1} 
when 
\begin{equation*}
|\lambda| \in (2,+\infty) ~~~\text{and}~~~\mu \in [-\lambda^2/4  , \lambda^2 / 4]
\end{equation*}
 or when
\begin{equation*}
|\lambda| \in (1,2]  ~~~~ \text{and}~~~\mu \in (1-|\lambda|,|\lambda|-1).
\end{equation*}

In 2018, X. Tang and W. Zhang (\cite{TX}) generalized the equation \eqref{general-1} to the one 
\begin{equation}
	\phi^{2}(x)=h(\phi(f(x)))+g(x),
	\label{equation}
\end{equation}
where $ h,f $ and $g $ are given functions, and $\phi $ is an unknown one.
By the Banach Contraction Principle, they gave the existence of bounded Lipschitz 
solutions on $\mathbb{R}$
of \eqref{equation} under Lipschitz condition in the case where $g$ is bounded and showed the existence of unbounded Lipschitz solutions on $\mathbb{R}$ of \eqref{equation} under additional bounded nonlinearities in the case where $g$ is unbounded. In addition, without Lipschitz condition, they applied piecewise construction method to give continuous solutions on $\mathbb{R}$.

In this paper, we proceed to study $C^1$ solutions   of the equation \eqref{equation}. In section 2, imposing conditions on derivatives of given functions, by the  Fiber Contraction Theorem, we prove the existence of $C^1$ solutions with bounded derivatives. Section 3 is devoted to the proof of some details related to the Fiber Contraction of Theorem. We also give an example to demonstrate our result in section 4.

We state the Fiber Contraction Theorem in the following for our convenience, which can be found in \cite {book}. 
\begin{lm}[Fiber Contraction Theorem]
Let $X$ and $Y$ be complete metric spaces. Assume that $\Gamma:X \times Y \rightarrow X \times Y $ defined by 
$$
\Gamma (x,y) = (\Lambda(x), \Phi(x,y) )
$$
is continuous, where 
$\Lambda: X\to X$ is contractive and $\Phi:X\times Y\to Y$ is uniformly contractive with respect to the first variable, that is,
$$
d_Y(\Phi(x,y_1), \Phi(x,y_2)) \le \gamma d_Y(y_1,y_2),
$$
in which $d_Y(\cdot, \cdot)$ denotes the metric in $Y$ and $0<\gamma<1$ is a constant. Then $\Lambda$ has a globally attracting fixed point $x_\infty$ in $X$ and $\Phi(x_\infty,\cdot)$ also has a globally attracting fixed point $y_\infty$ in $Y$.
In addition, $(x_\infty,y_\infty)$ is a globally attracting fixed point of $\Gamma$.
\label{lemma2}
\end{lm}

\section{Main results}
In this section, we give a result concerning the existence of $C^1$ solutions of the equation \eqref{equation}, which is stated in the following theorem.

\begin{thm}
Assume that functions $h:\mathbb{R}\rightarrow\mathbb{R}$, $f: \mathbb{R}\rightarrow\mathbb{R}$  and  $g:\mathbb{R}\rightarrow\mathbb{R}$ are of class $C^1$ such that
\begin{align}
\inf_{x\in \mathbb{R}}|h'(x)|\ge K,~~~~~~& \inf_{x\in \mathbb{R}}|f'(x)|\ge \alpha,
\label{derivatives-hf}
\\
\sup_{x\in\mathbb{R}}|g(x)|<+\infty,~~~~& \sup_{x\in\mathbb{R}}|g'(x)|\le \beta,
\label{condition-g}
\end{align}
where $K>1$, $\alpha>0$ and $\beta>0$ are given real constants such that
\begin{align}
\beta <\frac{1}{4}\alpha^2K^2  ~~~~~&\text{when}~~~ \alpha< 2(1-\frac{1}{K}),
\label{condi-cons1}
\\
\beta<(K-1)(\alpha K-K+1)~~~~~&\text{when}~~~ \alpha\ge 2(1-\frac{1}{K}).
\label{condi-cons2}
\end{align}
Then the functional equation \eqref{equation} has a  solution of class $C^1$ and its derivative is bounded.
\label{thm}
\end{thm}

\noindent
{\bf Proof.} Clearly,  it follows from the condition \eqref{derivatives-hf} that functions $h$ and $f$ are bijections on $\mathbb{R}$. Then the equation \eqref{equation} is equivalent to the form
\begin{equation}
	\phi(x)=h^{-1}(\phi^{2}(f^{-1}(x))-g(f^{-1}(x))), ~~~~~~x\in \mathbb{R}.
\end{equation}

Now we define some complete metric spaces.
Let $C^{0}_{b}(\mathbb{R}):=\{\phi:\mathbb{R}\rightarrow\mathbb{R}$ $\big|$ $\phi$ is continuous and $\sup_{x\in\mathbb{R}}|\phi(x)|<+\infty$\}. It is
evident that $C^{0}_{b}(\mathbb{R})$ is a Banach space equipped with the supremum norm ($\|\phi\|:=\sup_{x\in\mathbb{R}}|\phi(x)|$ for every $\varphi\in C^{0}_{b}(\mathbb{R})$).
For a constant $L\ge 0$, let
$C^{0}_{b}(\mathbb{R};L):=C^{0}_{b}(\mathbb{R})\cap\{\phi:\mathbb{R}\rightarrow\mathbb{R}|\text{Lip}(\phi)\leq L\}$, which is a nonempty closed subset. Thus, $C^{0}_{b}(\mathbb{R};L)$ is a complete metric space with the supremum norm in $C^{0}_{b}(\mathbb{R})$. 
In addition, for a constant $\rho>0$, we consider the set $\mathbb{F}_{\rho}:=C^{0}_{b}(\mathbb{R}) \cap \{\phi:\mathbb{R}\to \mathbb{R}|\|\phi\|\le \rho \}$, also a nonempty closed subset of $C^{0}_{b}(\mathbb{R})$ and a complete metric space with the supremum norm in $C^{0}_{b}(\mathbb{R})$.

We define  a bundle map
\begin{equation*}
	\Gamma:C^{0}_{b}(\mathbb{R},L)\times\mathbb{F_{\rho}}\rightarrow C^{0}_{b}(\mathbb{R},L)\times\mathbb{F_{\rho}}
\end{equation*}
by
\begin{equation}
	\Gamma(\phi,\Phi)=(\Lambda(\alpha),\Psi(\phi,\Phi)) ~~~\text{for}~~~(\phi,\Phi)\in C^{0}_{b}(\mathbb{R},L)\times\mathbb{F_{\rho}},
\end{equation}
where $\Lambda: C^{0}_{b}(\mathbb{R},L) \to C^{0}_{b}(\mathbb{R},L)$ is defined by 
\begin{equation}
	\Lambda(\phi)=h^{-1}\circ(\phi^{2}\circ f^{-1}-g\circ f^{-1}) ~~~\text{for all} ~~\phi\in  C^{0}_{b}(\mathbb{R},L),
	\label{Lambda}
\end{equation}
and $\Psi: C^{0}_{b}(\mathbb{R},L)\times\mathbb{F_{\rho}}\to \mathbb{F_{\rho}}$ is defined by
\begin{align}
\Psi(\phi,\Phi)&=(h^{-1})'\circ(\phi^2\circ f^{-1}-g\circ f^{-1})\cdot
\nonumber\\
&
~\quad\{\Phi\circ\phi\circ f^{-1} \cdot\Phi\circ f^{-1}-g'\circ f^{-1} \}\cdot (f^{-1})'
\label{Psi}
\end{align}
for all $(\phi,\Phi)\in
C^{0}_{b}(\mathbb{R},L)\times\mathbb{F_{\rho}}$, where $\cdot$ denotes multiplication of two functions, i.e., $(\phi_1\cdot\phi_2)(x):=\phi_1(x)\cdot\phi_2(x)$ for all $x\in\mathbb{R}$.

In what follows, $L$ is chosen to satisfy
\begin{align}
		\frac{1}{2}\alpha K-\frac{1}{2}\sqrt{\alpha^2K^2-4\beta}\leq &L\leq \frac{1}{2}\alpha K+\frac{1}{2}\sqrt{\alpha^2K^2-4\beta},
	\label{selfmapping}
	\\
	& L< K-1,
	\label{l-contraction}
\end{align}
and $\rho$ is chosen to satisfy 
\begin{align}
\frac{1}{2}\alpha K-\frac{1}{2}\sqrt{\alpha^2K^2-4\beta}\leq& \rho\leq \frac{1}{2}\alpha K+\frac{1}{2}\sqrt{\alpha^2K^2-4\beta},
\label{psi-selfmapiing}
\\
&\rho<\frac{1}{2}\alpha K.
\label{psi-contraction}
\end{align}
By \eqref{condi-cons1} or \eqref{condi-cons2}, such $L$ and $\rho$ exist. In fact, in the case \eqref{condi-cons1}, it is easy to check that $\alpha< 2(1-\frac{1}{K})$ implies
\begin{equation}
\frac{1}{2}\alpha K-\frac{1}{2}\sqrt{\alpha^2K^2-4\beta}<K-1.
\label{intersection-L}
\end{equation}
It follows that 
$$
\bigg[\frac{1}{2}\alpha K-\frac{1}{2}\sqrt{\alpha^2K^2-4\beta}, \frac{1}{2}\alpha K+\frac{1}{2}\sqrt{\alpha^2K^2-4\beta}\bigg] \cap (0,K-1)\ne \emptyset,
$$
which yields that there exists the $L$ satisfying \eqref{selfmapping} and \eqref{l-contraction}. Clearly,
\begin{equation}
\frac{1}{2}\alpha K-\frac{1}{2}\sqrt{\alpha^2K^2-4\beta}<\frac{1}{2}\alpha K
\label{intersection-rho}
\end{equation}
if $\beta<\frac{1}{4}\alpha^2K^2$. As a result, the $\rho$ satisfying \eqref{psi-selfmapiing} and \eqref{psi-contraction} exists. Note that $\beta<\frac{1}{4}\alpha^2K^2$ is required only to guarantee $\sqrt{\alpha^2K^2-4\beta}$ is positive.
In the other case \eqref{condi-cons2}, we see that 
$$
\beta<(K-1)(\alpha K-K+1) \le \bigg( \frac{K-1+\alpha K-K+1}{2}\bigg)^2 = \frac{1}{4}\alpha^2K^2,
$$
from which we also obtain \eqref{intersection-rho}. Therefore, the $\rho$ satisfying \eqref{psi-selfmapiing} and \eqref{psi-contraction} also exists. Moreover, when $\alpha\ge 2(1-\frac{1}{K})$, it is easy to calculate that $\beta<(K-1)(\alpha K-K+1)$ is equivalent to \eqref{intersection-L}. Consequently, the $L$ satisfying \eqref{selfmapping} and \eqref{l-contraction} can be chosen.

We claim that maps $\Lambda$ and $\Psi$ are well defined by \eqref{Lambda} and \eqref{Psi} under \eqref{selfmapping}-\eqref{psi-contraction}.
In fact, since $\phi$, $h^{-1}$,  $f^{-1}$  and  $g$ are all continuous on $\mathbb{R}$, so is $\Lambda(\phi)$ for each $\phi\in C^{0}_{b}(\mathbb{R},L)$.
Letting $M_{*}:=\max\{{\|\varphi\|},{\|g\|}\}$, it follows that
\begin{equation*}\begin{split}
\sup_{x\in\mathbb{R}}|\Lambda(\phi)(x)|&=\sup\limits_{x\in\mathbb{R}}|h^{-1}(\phi^{2}(f^{-1}(x))-g(f^{-1}(x)))|
\\
 &=\sup\limits_{x\in\mathbb{R}}|h^{-1}(\phi^{2}(x)-g(x))|\\
 &\le\sup\limits_{|x|\leq 2M_{*}}|h^{-1}(x)|<+\infty,	
	\end{split}
\end{equation*}
that is, $\Lambda(\phi)$ is also bounded on $\mathbb{R}$. Using the mean value theorem, by the first inequality of \eqref{derivatives-hf}, we derive that 
$$|h^{-1}(x)-h^{-1}(y)| = |(h^{-1})'(\xi)|\cdot |x-y| = \frac{1}{h'(\xi) } |x-y|\le \frac{1}{K}|x-y|$$
for all $x,y\in \mathbb{R}$, where $\xi$ is a point between $x$ and $y$. Similarly, by the first inequality of \eqref{derivatives-hf} and the second one of \eqref{condition-g}, we deduce respectively that 
\begin{equation*}
|f^{-1}(x)-f^{-1}(y)|\le \frac{1}{\alpha}|x-y|~~~\text{for all}~~x,y\in\mathbb{R}
\end{equation*}
and
\begin{equation*}
	|g(x)-g(y)|\le \beta|x-y|~~~\text{for all}~~x,y\in\mathbb{R}.
\end{equation*}
Therefore, for any $x_{1},x_{2}\in \mathbb{R}$, by \eqref{selfmapping}, we have
\begin{equation*}
	\begin{split}
		&|\Lambda(\phi)(x_{1})-\Lambda(\phi)(x_{2})|
		\\
		&=|h^{-1}(\phi^{2}(f^{-1}(x_1))-g(f^{-1}(x_1)))-h^{-1}(\phi^{2}(f^{-1}(x_2))-g(f^{-1}(x_2)))|\\
		&\leq\frac{1}{K}|\phi^{2}(f^{-1}(x_1))-g(f^{-1}(x_1)))-\phi^{2}(f^{-1}(x_2)+g(f^{-1}(x_2)))|\\
		&\leq\frac{1}{K}(\frac{L^2}{\alpha}+\frac{\beta}{\alpha})|x_1-x_2|\\
		&\leq L|x_1-x_2|,
	\end{split}
\end{equation*}
which shows Lip$(\Lambda(\phi))\le L$ for each $\phi\in C^{0}_{b}(\mathbb{R},L)$. It is shown that $\Lambda$ maps $C^{0}_{b}(\mathbb{R},L)$ into itself if \eqref{selfmapping} holds, that is, $\Lambda$ is well defined under \eqref{selfmapping}.
Evidently, $\Psi(\phi,\Phi)$ is continuous on $\mathbb{R}$ for every $(\phi,\Phi)\in C^{0}_{b}(\mathbb{R},L)\times\mathbb{F_{\rho}}$. 
Moreover, by \eqref{derivatives-hf}, \eqref{condition-g} and \eqref{psi-selfmapiing},
\begin{align*}
\|\Psi(\phi,\Phi)\| &= \sup_{x\in\mathbb{R}}|(h^{-1})'(\phi^2\circ f^{-1}(x)-g\circ f^{-1}(x))\cdot
\\
&\quad\{\Phi\circ\varphi\circ f^{-1}(x)\cdot\Phi\circ f^{-1}(x)-g'\circ f^{-1}(x)\} \cdot (f^{-1})'(x)|\\
&\le 
\sup_{x\in\mathbb{R}}|(h^{-1})'(x)| \bigg((\sup_{x\in\mathbb{R}}|\Phi(x)|)^2 \cdot \sup_{x\in\mathbb{R}}|(f^{-1})'(x)| 
\\
&\quad+ \sup_{x\in\mathbb{R}}|g'(x)|\cdot \sup_{x\in\mathbb{R}}|(f^{-1})'(x)|\bigg)
\\
&\le
\frac{\rho^2+\beta}{K\alpha}\le\rho,
\end{align*}
which shows $\Phi$ is well defined from $C^{0}_{b}(\mathbb{R},L)\times\mathbb{F_{\rho}}$ to $\mathbb{F_{\rho}}$.

We further have the following three assertions, whose proof will be given in next sections.
\begin{description}
\item [(A1)] $\Lambda$ is a contraction on $C^{0}_{b}(\mathbb{R},L)$;

\item [(A2)] $\Psi$ is a uniform contraction with respect to the first variable;

\item [(A3)] $\Gamma$ is continuous.
\end{description}
Therefore, by Lemma \ref{lemma2}, $\Gamma$ has a globally attracting fixed point $(\phi_*,\Phi_*)$, that is, for every $(\phi,\Phi)\in C^{0}_{b}(\mathbb{R},L)\times\mathbb{F_{\rho}}$, $\Gamma^n(\phi,\Phi)$ converges to $(\phi_*,\Phi_*)$ as $n\to +\infty$.
Choose arbitrarily $\phi_0\in C^{0}_{b}(\mathbb{R},L)$ and $\Phi_0\in \mathbb{F_{\rho}}$ such that $\Phi_0=(\phi_0)'$. Let
$$
(\phi_n,\Phi_n):= \Gamma^n(\phi_0,\Phi_0).
$$
In accordance with the definitions of $\Lambda$ and $\Psi$, we have that 
$(\phi_n)'=\Phi_n$ for all $n\ge 0$.
Furthermore, $(\phi_n,\Phi_n)\to (\phi_*,\Phi_*)$ as $n\to +\infty$.
It follows that $(\phi_*)'=\Phi_*$, implying that $\phi_*$ is of class $C^1$ and its derivative is bounded by $\rho$. Note that $\phi_*$ is a fixed point of $\Lambda$, which is a solution of the equation \eqref{equation}.
The proof is completed.
\hfill
$\square$

\section{Proofs of assertions (A1)-(A3)}

{\bf Proof of assertion (A1)}. 
If $\phi_1$ and $\phi_2$ are in the set $C^{0}_{b}(\mathbb{R},L)$, by the first inequality of \eqref{derivatives-hf},
\begin{equation*} \begin{split}
&\|\Lambda(\phi_1)-\Lambda(\phi_2)\|
  =\sup\limits_{x\in\mathbb{R}}|\Lambda\phi_1(x)-\Lambda\phi_2(x)|\\
  &=\sup\limits_{x\in\mathbb{R}}|h^{-1}(\phi_1^2(f^{-1}(x))-g(f^{-1}(x)))-h^{-1}(\phi_2^2(f^{-1}(x))-g(f^{-1}(x)))|\\
  &\leq\frac{1}{K}\sup\limits_{x\in\mathbb{R}}|\phi_1^2\circ f^{-1}(x)-\phi_2^2\circ f^{-1}(x)|\\
  &=\frac{1}{K}\sup\limits_{x\in\mathbb{R}}|\phi_1^2(x)-\phi_2^2(x)|\\
  &\leq\frac{1}{K}{\sup\limits_{x\in\mathbb{R}}|\phi_1^2(x)-\phi_1(\phi_2(x))|+\sup\limits_{x\in\mathbb{R}}|\phi_1(\phi_2(x))-\phi_2^2(x)|}\\
  &\leq\frac{1}{K}(L+1)\|{\phi_1-\phi_2}\|.
\end{split} 
\end{equation*}   
As a result, it follows from \eqref{l-contraction} that $\Lambda$ is a contraction on $C^{0}_{b}(\mathbb{R},L)$. The proof of assertion {\bf (A1)} is completed.
\hfill
$\square$

\noindent
{\bf Proof of assertion (A2).} 
By \eqref{derivatives-hf}, for any $\phi\in C^{0}_{b}(\mathbb{R},L)$ and any  $\Phi_1,\Phi_2\in\mathbb{F_{\rho}}$, we have that
\begin{equation*}\begin{split}
		&\|\Psi(\phi,\Phi_{1})-\Psi(\phi,\Phi_{2})\|
		=\sup_{x\in\mathbb{R}}|\Psi(\phi,\Phi_{1})(x)-\Psi(\phi,\Phi_{2})(x)|\\
		&=\sup_{x\in\mathbb{R}}|(h^{-1})'(\phi^2\circ f^{-1}(x)-g\circ f^{-1}(x))\cdot\{\Phi_{1}\circ\phi\circ f^{-1}(x)\cdot\\
		&\quad \Phi_{1}\circ f^{-1}(x)-\Phi_{2}\circ\phi\circ f^{-1}(x)\cdot\Phi_{2}\circ f^{-1}(x)\} \cdot(f^{-1})'(x)|\\
		&\leq\frac{1}{K}\cdot\frac{1}{\alpha}\cdot \sup_{x\in\mathbb{R}}|\Phi_{1}\circ\phi\circ f^{-1}(x)\cdot
		\Phi_{1}\circ f^{-1}(x)-\Phi_{2}\circ\phi\circ f^{-1}(x)\cdot\Phi_{2}\circ f^{-1}(x)|\\
		&\leq\frac{1}{\alpha K} \sup_{x\in\mathbb{R}}|\Phi_{1}\circ\phi\circ f^{-1}(x)\cdot
		\Phi_{1}\circ f^{-1}(x)-\Phi_{1}\circ\phi\circ f^{-1}(x)\cdot\Phi_{2}\circ f^{-1}(x)|\\
		&\quad+\frac{1}{\alpha K} \sup_{x\in\mathbb{R}}|\Phi_{1}\circ\phi\circ f^{-1}(x)\cdot
		\Phi_{2}\circ f^{-1}(x)-\Phi_{2}\circ\phi\circ f^{-1}(x)\cdot\Phi_{2}\circ f^{-1}(x)|\\
		&\leq\frac{2\rho}{\alpha K}\|\Phi_{1}-\Phi_{2}\|.
	\end{split} 
\end{equation*} 
Therefore, $\Psi$ is a uniform contraction by \eqref{psi-contraction}.
The proof of assertion {\bf (A2)} is completed.
\hfill
$\square$

\noindent
{\bf Proof of assertion (A3).} 
To prove the continuity of $\Gamma$, it suffices to show that the function $\phi\mapsto\Psi(\phi,\Phi_{0})$ is continuous for any fixed $\Phi_0\in \mathbb{F_{\rho}}$ since $\Lambda$ is a contraction and $\Psi$ is a uniform contraction with respect to the second variable $\Phi$.

Let
\begin{align*}
	{\cal S}(\phi)&:=(h^{-1})'\circ (\phi^2\circ f^{-1}-g\circ f^{-1}),
	\\
	{\cal W}(\phi,\Phi)&:=\Phi\circ\phi\circ f^{-1}\cdot\Phi\circ f^{-1}\cdot (f^{-1})',
	\\
	{\cal C}&:= g'\circ f^{-1}\cdot (f^{-1})'.
\end{align*}
Then by \eqref{Psi} $\Psi(\phi,\Phi)$ can be rewritten as
\begin{equation}
	\Psi(\phi,\Phi)={\cal S}(\phi)\cdot {\cal W}(\phi,\Phi)-{\cal S}(\phi)\cdot {\cal C}.
\end{equation}
It is easy to calculate that
\begin{align}
	&\|\Psi(\phi,\Phi_{0})-\Psi(\phi_{0},\Phi_{0})\|
	\nonumber\\
	&=\|{\cal S}(\phi)\cdot {\cal W}(\phi,\Phi_{0})-{\cal S}(\phi)\cdot {\cal C}-{\cal S}(\phi_{0})\cdot {\cal W}(\phi_{0},\Phi_{0})+{\cal S}(\phi_{0})\cdot {\cal C}\|
	\nonumber\\
	&\leq\|{\cal S}(\phi)\cdot {\cal W}(\phi,\Phi_{0})-{\cal S}(\phi_{0})\cdot {\cal W}(\phi_{0},\Phi_{0})\|+ \|{\cal S}(\phi)\cdot {\cal C}-{\cal S}(\phi_{0})\cdot {\cal C}\|
	\nonumber\\
	&\le\|{\cal S}(\phi)\cdot {\cal W}(\phi,\Phi_{0})-{\cal S}(\phi)\cdot {\cal W}(\phi_{0},\Phi_{0})\|+\|{\cal S}(\phi)\cdot {\cal W}(\phi_{0},\Phi_{0})-  {\cal S}(\phi_{0})\cdot {\cal W}(\phi_{0},\Phi_{0})\|
	\nonumber\\
	&\quad + \| {\cal S}(\phi)-{\cal S}(\phi_0)\|\cdot \|{\cal C}\|
	\nonumber
\\
&\leq \|{\cal S}(\phi)\|\cdot \| {\cal W}(\phi,\Phi_{0})-{\cal W}(\phi_{0},\Phi_{0})\|+(\|{\cal W}(\phi_{0},\Phi_{0})\|+\|{\cal C}\|)\cdot\| {\cal S}(\phi)-{\cal S}(\phi_0)\|
\nonumber\\
&\le \frac{1}{K}  \| {\cal W}(\phi,\Phi_{0})-{\cal W}(\phi_{0},\Phi_{0})\| + \bigg(\frac{\rho^2+\beta}{\alpha}\bigg)\| {\cal S}(\phi)-{\cal S}(\phi_0)\|.
\label{calculate}
\end{align}
As a consequence, it is sufficient to show that ${\cal W}(\phi,\Phi_0)$ and ${\cal S}(\phi)$ both are continuous at $\phi_0$. We need the following lemma.

\begin{lm}
Maps 
	\begin{align*}
		\phi\in C^0_b(\mathbb{R})\mapsto {\cal T}(\phi)&:=\Phi_0\circ \phi \circ f^{-1}\in C^0_b(\mathbb{R})~~~\text{and}\\
		\phi\in C^0_b(\mathbb{R},L)\mapsto {\cal S}(\phi)&= (h^{-1})'\circ (\phi^2 \circ f^{-1} - g\circ f^{-1})\in C^0_b(\mathbb{R})
	\end{align*}
	are continuous.
	\label{lemma}
\end{lm}
The proof of Lemma \ref{lemma} is given after we complete the proof of assertion {\bf (A3)}. It is easy to see that 
$$
\| {\cal W}(\phi,\Phi_{0})-{\cal W}(\phi_{0},\Phi_{0})\| \le \frac{\rho}{\alpha} \|\Phi_0\circ\phi\circ f^{-1}-\Phi_0\circ\phi_{0}\circ f^{-1}\|.
$$
Thus, by the continuity of ${\cal T}$ in Lemma \ref{lemma}, we obtain that ${\cal W}(\phi,\Phi_{0})$ is continuous at $\phi_{0}$. Consequently, by \eqref{calculate} and the continuity of ${\cal W}(\phi,\Phi_{0})$ and ${\cal S}(\phi)$, $\Psi(\phi,\Phi_{0})$ is continuous at $\phi_0$. The proof is completed.
\hfill
$\square$

\noindent
{\bf Proof of Lemma \ref{lemma}.} 
We first prove that ${\cal T}$ is continuous.
Fix $\phi_{0}\in C^{0}_{b}(\mathbb{R})$ arbitrarily and we need to show that ${\cal T}$ is continuous at $\phi_0$. In other words, for every $\epsilon>0$ there exists a $\delta>0$ such that
\begin{equation*}
\|\Phi_0\circ\phi\circ f^{-1}-\Phi_0\circ\phi_{0}\circ f^{-1}\| < \epsilon~~ whenever~~ \|\phi-\phi_{0}\| < \delta.
\end{equation*}
Since the function $\Phi_0$ is continuous on $\mathbb{R}$, it is uniformly continuous on the bounded closed interval $I_{{\cal T}}:=[-\|\phi_0\|-1,\|\phi_0\|+1]$.
Namely, for every $\epsilon>0$ there exists a $0<\delta_0<1$ such that 
$|\Phi_0(x_1)-\Phi_0(x_2)|<\epsilon/2$ whenever $x_1,x_2\in I_{{\cal T}}$ and $|x_1-x_2|<\delta_0$.
Notice that, when $\|\phi-\phi_0\|<\delta_0$,
$$|\phi\circ f^{-1}(x)| \le \sup_{x\in \mathbb{R}}|\phi(x)|=\|\phi\|\le \|\phi-\phi_0\|+\|\phi_0\|<1+\|\phi_0\|$$ for all $x\in \mathbb{R}$, that is, $\phi\circ f^{-1}(x)\in I_{{\cal T}}$ for all $x\in \mathbb{R}$, and 
$$
|\phi\circ f^{-1}(x)-\phi_0\circ f^{-1}(x)| \le \sup_{x\in \mathbb{R}} |\phi\circ f^{-1}(x)-\phi_0\circ f^{-1}(x)| \le \|\phi-\phi_0\|<\delta_0
$$
for all $x\in \mathbb{R}$. In particular, $\phi_0\circ f^{-1}(x)\in I_{{\cal T}}$ for all $x\in \mathbb{R}$. Thus, whenever $\|\phi-\phi_0\|<\delta_0$, 
\begin{equation*}
|\Phi_0\circ \phi\circ f^{-1}(x)-\Phi_0\circ\phi_0\circ f^{-1}(x)|<\epsilon/2~~~\text{for all}~~x\in\mathbb{R},
\end{equation*}
which implies that $\|\Phi_0\circ\phi\circ f^{-1}-\Phi_0\circ\phi_{0}\circ f^{-1}\|= \sup_{x\in\mathbb{R}} |\Phi_0\circ \phi\circ f^{-1}(x)-\Phi_0\circ\phi_0\circ f^{-1}(x)| \le \epsilon/2<\epsilon$. The continuity of ${\cal T}$ is proved.

Next, we prove that $\cal S$ is also continuous. Fixing arbitrarily $\phi_0\in C^{0}_{b}(\mathbb{R})$, we need to show that $\cal S$ is continuous at $\phi_0$. Namely, for every $\epsilon>0$, there exists a $\delta>0$ such that
\begin{equation*}
\|(h^{-1})'\circ(\phi^2\circ f^{-1}-g\circ f^{-1})-(h^{-1})'\circ(\phi_0^2\circ f^{-1}-g\circ f^{-1})\| < \epsilon~~ whenever~~ \|\phi-\phi_{0}\| < \delta.
\end{equation*}
Since the function $(h^{-1})'$ is continuous on $\mathbb{R}$, it is uniformly continuous on the bounded closed interval $I_{{\cal S}}:=[-\|\phi_0\|-\|g\|-1,\|\phi_0\|+\|g\|+1]$. In other words, for every $\epsilon>0$ there exists a $0<\delta_0<1$ such that 
$|(h^{-1})'(x_1)-(h^{-1})'(x_2)|<\epsilon/2$ whenever $|x_1-x_2|<\delta_0$ and $x_1,x_2\in I_{{\cal S}}$.
When $\|\phi-\phi_0\|<\delta_0/(L+1)$,
\begin{align*}
|\phi^2\circ f^{-1}(x)-g\circ f^{-1}(x)| &\leq \sup_{x\in \mathbb{R}}|\phi^2\circ f^{-1}(x)-g\circ f^{-1}(x)|
\\
&=\|\phi^2\circ f^{-1}-g\circ f^{-1}\|\leq \|\phi\|+\|g\|
\\
&=\|\phi-\phi_0\|+\|\phi_0\|+\|g\|
\\
& \leq 1+\|\phi_0\|+\|g\|
\end{align*}
for all $x\in \mathbb{R}$, that is, $\phi^2\circ f^{-1}(x)-g\circ f^{-1}(x)\in I_{\cal S}$
for all $x\in \mathbb{R}$, and 
\begin{align*}
|\phi^2\circ f^{-1}(x)-\phi_0^2\circ f^{-1}(x)| &\le \| \phi^2 - \phi_0^2\| \le \| \phi^2 -\phi\circ \phi_0\| +\|\phi\circ \phi_0-\phi_0^2\|
\\
&\le (L+1)\|\phi-\phi_0\|<\delta_0
\end{align*}
for all $x\in \mathbb{R}$. Thus, whenever $\|\phi-\phi_0\|<\delta_0/(L+1)$, 
\begin{equation*}
|(h^{-1})'\circ(\phi^2\circ f^{-1}-g\circ f^{-1})(x)-(h^{-1})'\circ(\phi_0^2\circ f^{-1}-g\circ f^{-1})(x)|<\frac{\epsilon}{2}
\end{equation*}
for all $x\in \mathbb{R}$, which implies that $\|(h^{-1})'\circ(\phi^2\circ f^{-1}-g\circ f^{-1})-(h^{-1})'\circ(\phi_0^2\circ f^{-1}-g\circ f^{-1})\| \le \epsilon/2<\epsilon$. The continuity of ${\cal S}$ is proved.
\hfill
$\square$

\section{An example}

In this section we give an example to demonstrate our Theorem \ref{thm}.

{\bf Example:} Our Theorem \ref{thm} can be applied to the equation
\begin{equation}
	\phi^2(x)= \sin(\phi(e^x+5x))+4\phi(e^x+5x)+\cos x,
	\label{example}
\end{equation}
which is of the form \eqref{equation} with $f(x)=e^x+5x$, $h(x)=\sin x+4x$ and $g(x)=\cos x$. One can check that $f,g$ and $h$ satisfy the conditions $\eqref{derivatives-hf}$ and $\eqref{condition-g}$ with constants $\ K=3$, $\alpha=5 $, $\beta=1$. Furthermore, $2(1-1/\ K)=1< \alpha$ and
\begin{equation*}
	(K-1)(\alpha K-K+1)=26>\beta,
\end{equation*}
i.e., condition \eqref{condi-cons2} is fulfilled. Consequently, by Theorem \ref{thm}, the equation \eqref{example} has a  solution of class $C^1$ and its derivative is bounded.



\begin{thebibliography}{99}

\bibitem{Balcerowski} 
M. Balcerowski, On the functional equation $x+f(y+f(x)) = y+f(x+f(y))$,
{\it Aequationes Math.} {\bf 75} (2008), 297-303.

\bibitem{Baron}    
K. Baron, Recent results in the theory of functional equations in a single variable, Survey in 40th ISFE (August 2002, 
Gron\'{o}w, Poland); Seminar LV, No. 15, 
{\it Mathematisches Institut \uppercase\expandafter{\romannumeral1}, Universit\"{a}t Karlsruhe,} 2003,
http://www.mathematik.uni-karlsruhe.de$/^\sim$ semlv/.

\bibitem{BJ}
K. Baron and  W. Jarczyk, Recent results on functional equations in a single variable,
{\it Aequat. Math.} {\bf 61} (2001), 1-48.


\bibitem{Belluot}  
N. Brillou\"{e}t-Belluot, Problem 15, Proceedings of 38th ISFE (2000 Hungary), 
{\it Aequationes Math.} {\bf 61} (2001), 304.


\bibitem{ZWN}     
N. Brillou\"et-Belluot and W. Zhang, On a class of iterative-difference equations, 
{\it J. Difference Eq. Appl.} {\bf 16} (2010), no. 11, 1237-1255.

\bibitem{book}
C. Chicone,  {\it Ordinary Differential Equations with Applications}, Springer, New York, 1999.


\bibitem{DM}
S. Draga and J. Morawiec, Reducing the polynomial-like iterative equation order and a generalized Zolt\'an Boros' problem, {\it Aequationes Math.} {\bf 90} (2016), 935-950.

\bibitem{Yarczyk}  
J. Jarczyk and W. Jarczyk, On a problem of N. Brillou\"{e}t-Belluot, 
{\it Aequationes Math.} {\bf 72} (2006), 198-200.

\bibitem{KCG}
M. Kuczma,  B. Choczewski and R. Ger, {\it Iterative Functional Equations}, Cambridge University
Press, Cambridge,  1990.

\bibitem{LJLZ}
L. Liu, W. Jarczyk, L. Li and W. Zhang, Iterative roots of piecewise monotonic functions
of nonmonotonicity height not less than 2, {\it  Nonlinear Anal.} {\bf 75} (2012), 286-303.

\bibitem{LZ}
L. Liu and W. Zhang,Non-monotonic iterative roots extended from characteristic intervals,
{\it J. Math. Anal. Appl.} {\bf 378} (2011), 359-373. 

\bibitem{Matkowski}      
J. Matkowski and W. Zhang, 
Method of characteristic for functional equations in polynomial form, 
{\it Acta Math. Sinica, New Series} {\bf 13} (1997), no. 3, 421-432.

\bibitem{MR}
A. Mukherjea and J. S. Ratti, On a functional equation involving iterates of a bijection on
the unit interval, {\it Nonlinear Anal.} {\bf 7} (1983), 899-908.

\bibitem{Sablik}   
M. Sablik, Remark 4, Report of meeting. 10th International Conference on Functional Equations and Inequalities (B\c{e}dlewo, 2005), 
{\it Ann. Acad. Paed. Cracoviensis Studia Math.} {\bf 5} (2006), 
158.

\bibitem{TX}      
X. Tang and WN. Zhang,  Continuous solutions of a second order iterative equation,
{\it Publ. Math. Debrecen} {\bf 93} (2018), 303-321.

\bibitem{YDL}        
D. Yang and W. Zhang, Characteristic solutions of polynomial-like iterative equations, 
{\it Aequationes Math.} {\bf 67} (2004), 80-105.

\bibitem{ZYY}      
Y. Zeng and W. Zhang,  Continuous solutions of an iterative-difference equation and Brillou\"{e}t-Belluot's problem,
{\it Publ. Math. Debrecen} {\bf 78} (2011), no. 3-4, 613-624.


\bibitem{ZNX}
W. Zhang,  K. Nikodem and  B. Xu, Convex solutions of polynomial-like iterative equations,
{\it J. Math. Anal. Appl.} {\bf  315} (2006), 29-40.




\end{thebibliography}
\end{document}